\documentclass[11pt]{article}
\usepackage[margin=1in]{geometry} 
\UseRawInputEncoding
\usepackage{hyperref}
\usepackage[affil-sl]{authblk}
\usepackage{lineno}

\usepackage[utf8]{inputenc}
\usepackage[T1]{fontenc}
\usepackage{anyfontsize}
\usepackage{lipsum} 
\usepackage{microtype} 
\usepackage{titlesec} 
\usepackage{graphicx} 
\usepackage{svg}
\usepackage[normalem]{ulem}
\usepackage{marginnote} 

\usepackage{bbm}
\usepackage{xspace}
\usepackage{tabularx}
\usepackage{stmaryrd}
\SetSymbolFont{stmry}{bold}{U}{stmry}{m}{n}
\usepackage{bbm}
\usepackage{url}            
\usepackage{booktabs}       
\usepackage{amsfonts}       
\usepackage{nicefrac}       
\usepackage{microtype}      
\usepackage{amsmath,amsthm}
\usepackage{tikz}
\usepackage[most]{tcolorbox}
\usepackage{floatrow} 
\usepackage{nicematrix}
\usepackage[square,numbers]{natbib}

\usepackage{newtxmath}
\usepackage{letltxmacro}
\usepackage{bm}

\usepackage{soul}
\usepackage[linesnumbered,ruled,vlined]{algorithm2e}
\usepackage[compatible]{algpseudocode} 
\usepackage{bbm}
\usepackage{svg}
\usepackage{xspace}
\usepackage{stmaryrd}
\usepackage{bbm}
\usepackage{url}            
\usepackage{booktabs}       
\usepackage{nicefrac}       
\usepackage{microtype}      
\usepackage{letltxmacro}
\usepackage{bm}
\usepackage[font=small,labelfont={bf},textfont={tt}]{caption} 
\usepackage{subcaption} 
\usepackage{soul}
\usepackage[linesnumbered,ruled,vlined]{algorithm2e}
\usepackage[compatible]{algpseudocode} 
\usepackage{doi}
\usepackage{fancyvrb,cprotect}
\usepackage{centernot}
\usepackage{fancyhdr}
\graphicspath{ {assets/} }

\usepackage{hyperref} 

\hypersetup{
    colorlinks=true,
    linkcolor=kon-peki,
    citecolor=kon-peki, 
    urlcolor=tsutsuji 
}

\usepackage{cleveref}       

\usepackage{framed}
\usepackage{rotating}

\usepackage[colorinlistoftodos]{todonotes}

\UseRawInputEncoding
\usepackage{thmtools}
\usepackage[framemethod=TikZ]{mdframed}

\usepackage{marginnote} 
\usepackage{xcolor}     
\usepackage{lipsum}     

\definecolor{ama-iro}{RGB}{0, 158, 243.0} 
\definecolor{fuyu-gaki}{HTML}{C75146} 
\definecolor{momiji}{RGB}{245, 70, 111} 
\definecolor{hotaru-bi}{RGB}{229,221,58} 
\definecolor{kon-peki}{RGB}{1,120,217} 
\definecolor{shin-kai}{RGB}{77,98,152} 
\definecolor{shin-ryoku}{RGB}{1,145,97} 
\definecolor{yama-budo}{RGB}{171,14,122} 
\definecolor{tsutsuji}{HTML}{C71585} 
\definecolor{mizu}{rgb}{0.0, 0.6, 0.8} 
\definecolor{kikyou}{rgb}{0.4, 0.4, 0.8} 
\definecolor{theoremcolor}{RGB}{77,98,152}
\definecolor{definitioncolor}{RGB}{255, 213, 179}
\definecolor{examplecolor}{RGB}{255, 213, 179}
\definecolor{verylightgray}{gray}{0.88}

\newtheoremstyle{theoremstyle}
  {\topsep} 
  {\topsep} 
  {} 
  {} 
  {\bfseries\color{black}} 
  {} 
  {.5em} 
  {\thmname{#1}~\thmnumber{#2}:~\textcolor{theoremcolor}{\thmnote{[#3]}}} 

\mdfdefinestyle{theoremframe}{
  linecolor=theoremcolor!80,
  linewidth=2pt,
  backgroundcolor=theoremcolor!05,
  roundcorner=5pt,
  nobreak=true,
  topline=false,
  bottomline=false,
  rightline=false,
}
\theoremstyle{theoremstyle}
\newtheorem{theorem}{Theorem}[section]
\surroundwithmdframed[style=theoremframe]{theorem}

\newtheoremstyle{remarkstyle}
  {\topsep} 
  {\topsep} 
  {} 
  {} 
  {\bfseries\color{black}} 
  {} 
  {0.5em} 
  {} 

\mdfdefinestyle{remarkframe}{
  linecolor=theoremcolor!80,
  linewidth=2pt,
  backgroundcolor=theoremcolor!05,
  roundcorner=5pt,
  nobreak=true,
  topline=false,
  bottomline=false,
  rightline=false,
}
\theoremstyle{remarkstyle}

\surroundwithmdframed[style=remarkframe]{remark}

\newtheoremstyle{definitionstyle}
  {\topsep} 
  {\topsep} 
  {} 
  {} 
  {\bfseries\color{black}} 
  {} 
  {.5em} 
  {\thmname{#1}~\thmnumber{#2}:~\textcolor{theoremcolor}{\thmnote{[#3]}}} 

\mdfdefinestyle{definitionframe}{
  linecolor=definitioncolor!80,
  linewidth=2pt,
  backgroundcolor=definitioncolor!45,
  roundcorner=5pt,
  nobreak=true,
  topline=false,
  bottomline=false,
  rightline=false,
}
\theoremstyle{definitionstyle}
\newtheorem{definition}[theorem]{Definition}
\surroundwithmdframed[style=definitionframe]{definition}

\newtheoremstyle{lemmastyle}
  {\topsep} 
  {\topsep} 
  {} 
  {} 
  {\bfseries\color{black}} 
  {} 
  {.5em} 
  {\thmname{#1}~\thmnumber{#2}:~\textcolor{theoremcolor}{\thmnote{[#3]}}} 

\mdfdefinestyle{lemmaframe}{
  linecolor=theoremcolor!80,
  linewidth=2pt,
  backgroundcolor=theoremcolor!05,
  roundcorner=5pt,
  nobreak=true,
  topline=false,
  bottomline=false,
  rightline=false,
}
\theoremstyle{lemmastyle}
\newtheorem{lemma}[theorem]{Lemma}
\surroundwithmdframed[style=lemmaframe]{lemma}

\surroundwithmdframed[style=definitionframe]{claim}
\newtheorem{corollary}[theorem]{Corollary}
\surroundwithmdframed[style=remarkframe]{corollary}

\newtcolorbox[auto counter, number within=section]{problem}[2][]{colframe=fuyu-gaki, colback=mizu!10, coltitle=white, title={Problem }, sharp corners, boxrule=0.8mm, width=0.99\textwidth, boxsep=2mm, left=3mm, right=3mm, top=2mm, bottom=2mm, breakable}

\usepackage{accents}
\usepackage{calc}
\newcommand{\Def}{=}
\makeatletter
\newcommand{\smalldollar}{\mathrel{\mathpalette\small@dollar\relax}}
\newcommand{\small@dollar}[2]{%
  \vcenter{\hbox{%
    $#1\textnormal{\fontsize{0.7\dimexpr\f@size pt}{0}\selectfont\$}$%
  }}%
}
\makeatother
\renewcommand{\emph}[1]{\textit{#1}}
\newcommand\Bigger[2][7]{\left#2\rule{0mm}{#1truemm}\right.}
\renewcommand{\vec}[1]{\mathbf{#1}}

\renewcommand{\epsilon}{\varepsilon}

\newcommand{\RCloseLOpenInterval}[2]{}

\newcommand{\Interval}[2]{\left[\right]}
\newcommand{\OpenInterval}[2]{\left(\right)}

\newcommand{\Dist}{\mathcal{D}}

\newcommand{\Indicator}[1]{\mathbbm{1}\left\{#1\right\}}

\newcommand{\bit}{\{0,1\}}
\newcommand{\poly}{\mathsf{poly}}
\newcommand{\Naturals}{\mathbb{N}}

\newcommand{\Field}{\mathbb{F}}
\newcommand{\BigO}[1]{O\left(#1\right)}

\newcommand{\BigOmega}[1]{\Omega\left(#1\right)}

\newcommand{\Mean}[2]{\mathbb{E}_{#1}\left[#2\right]}
\newcommand{\Prob}[1]{\Pr\left[#1 \right]}
\newcommand{\PProb}[2]{\Pr_{#2}\left[#1 \right]}
\newcommand{\True}{\texttt{True}}
\newcommand{\False}{\texttt{False}}

\renewcommand{\P}{\mathsf{P}}
\newcommand{\NP}{\mathsf{NP}}
\newcommand{\CoNP}{\mathsf{coNP}}

\newcommand{\Proof}{\pi}
\newcommand{\Size}[1]{|#1|}
\newcommand{\Card}[1]{\text{Card}(#1)}
\newcommand{\Axioms}{\mathcal{Q}}
\newcommand{\axiom}{q}
\newcommand{\PM}[1]{\text{PM}(#1)}

\newcommand{\highlight}[1]{\textcolor{red}{#1}}
\newcommand{\red}[1]{#1}
\newcommand{\green}[1]{#1}

\newcommand{\GreenB}{\green{B}}
\newcommand{\RedA}{\red{A}}

\newcommand{\Graph}{G}
\newcommand{\Vertices}[1]{V(#1)}
\newcommand{\Edges}[1]{E(#1)}
\newcommand{\CutEdges}[3]{e_{#3}(#1, #2)}

\newcommand{\OddComponents}[1]{q(#1)}
\newcommand{\degree}[2]{\mathsf{deg}_{#1}(#2)} 
\newcommand{\MaxDegree}[1]{\Delta_{#1}}
\newcommand{\Neighbourhood}[2]{\Gamma_{#1}\left(#2\right)}

\newcommand{\HardInstance}{H}
\newcommand{\Embedding}[1]{\Psi\left(#1\right)}
\newcommand{\Path}[2]{\mathsf{P}_{#1 \rightsquigarrow #2
}}

\newcommand{\Degree}[1]{\text{Deg}\left(#1\right)}

\newcommand{\Complement}[1]{\overline{#1}}

\newcommand{\PC}{\vdash_{\texttt{PC}_F}}
\newcommand{\SOS}{\vdash_{\texttt{SOS}}}

\newcommand{\EnDeeLambda}{(n, d, \lambda)}
\newcommand{\Subdivision}[2]{{#1}^{#2}}
\newcommand{\MinimalDegree}[1]{\delta({#1})}

\newcommand{\ExpansionFactor}[1]{\lambda(#1)}
\newcommand{\EmbeddingFunc}{\Psi}
\newcommand{\cTilde}{\frac{9d}{10}}

\title{Refuting Perfect Matchings in Spectral Expanders is Hard}

\author[1]{Ari Biswas}
\author[2]{Rajko Nenadov\thanks{Research supported by the Marsden Fund of the Royal Society of New Zealand.}}
\affil[1]{\small University Of Warwick, United Kingdom}
\affil[2]{\small University Of Canterbury, New Zealand}

\date{}
\begin{document}

\maketitle
\begin{abstract}
This work studies the complexity of refuting the existence of a perfect matching in spectral expanders with an odd number of vertices, in the Polynomial Calculus (PC) and Sum of Squares (SoS) proof system.
Austrin and Risse [SODA, 2021] showed that refuting perfect matchings in sparse $d$-regular \emph{random} graphs, in the above proof systems, with high probability requires proofs with degree $\Omega(n/\log n)$. 
We extend their result by showing the same lower bound holds for \emph{all} $d$-regular graphs with a mild spectral gap.
\end{abstract}

\section{Introduction}


Perhaps the most fundamental problem in computation is to provide an answer to the question asked by Cook and Reckchow in their seminal paper \citep{cook1979relative} -- ``\textit{Given a true statement $A$, is there a short proof of the claim that $A$ is true?}''
In trying to answer this question, we must first describe what constitutes a valid proof. That is, we must describe the language in which the proof is written (axioms), and the rules for checking it (the verifier).
Each set of rules for writing and checking a proof defines a proof system.
Therefore, a precise restatement of the question above is the following:  ``\textit{Given a true statement $A$ and a proof system $S$, what is the length of a shortest $\Proof \in S$ that proves $A$?}''
If we could show that there exists a proof system $S$, such that for \emph{any} true statement $A$, the length of the shortest proof in $S$ is upper bounded by some polynomial in the length of $A$, it would imply that $\CoNP = \NP$ and consequently the polynomial hierarchy collapes to $\NP$. 
Conversely, if we could show large proof size lower bounds for some $A$ in \emph{all} proof systems, it would lead to a formal proof of the widely believed conjecture that $\P \neq \NP$.
Unfortunately, such lower  bounds for \emph{arbitrary} proof systems are out of reach.
As an intermediate step,
the research community has invested a significant effort in proving lower bounds for increasingly expressive proof systems (e.g.\ see \citep{abascal2021strongly, alekhnovich2001lower, atserias2020Size, buss1999linear,  conneryd2023graph, de2023clique, impagliazzo1999lower,  kothari2017SumOfSquares, potechin:LIPIcs.CCC.2020.38, raz2008elusive, razborov1998lower, schoenebeck2008linear}). \par
In this work, we focus on the algebraic and semi-algebraic proof systems\footnote{Similar to the work of Austrin and Risse \citep{Austrin_2022}, our lower bounds also immediately extend to bounded depth Frege proof systems. However, the key technical contribution in this work is the graph theoretic techniques proposed for embedding a hard instance into the host graph.
Provided we are able to find an appropriate hard instance (which is known to exist for Frege proof systems), our techniques are independent of the underlying proof systems.
Thus, we restrict our preliminaries to PC and SoS to keep our preliminaries stand-alone, and yet concise.} of \emph{polynomial calculus (PC)} and \emph{sum of squares (SoS)}. 
In algebraic proof systems, we are given a set $\Axioms =\{\axiom_i(\vec{x}) \text{ } |\text{ } i \in [m] \}$ of $m$ polynomial equations\footnote{Semi-algebraic proof systems also allow for inequalities but we will not deal with inequality constraints in this paper.} over $n$ variables $\vec{x} = \{x_1, \dots, x_n\}$.
In PC, the equations can be over an arbitrary, but fixed field $\Field$, and in the SoS the coefficients are  over the reals.
We say a proof $\Proof$ is a refutation of $\Axioms$, if it is a proof of the claim (in the specified language) that there exists no assignment of $\vec{x} \in \Field^n$ that satisfies \emph{all} the polynomial equations in $\Axioms$.
In PC and SoS, the proof $\Proof$ is itself expressed as a sequence of polynomials.
Two common measures of the complexity of a semi-algebraic proof are size (the number of monomials appearing in the proof) and the degree (the largest degree of the proof polynomials that refute $Q$).
Trade-offs between the two are known.
If $\Axioms$ is the indexed set of polynomials with degree at most $k$, and if $d$ is the degree of the SoS refutation; then Atserias and Hakoniemi \cite[Corollary 1]{atserias2020Size} show that the size $s$ of the refutation is related as $s \ge \exp\left(\frac{(d-k-4)^2}{32(n+1)}\right)$. 
In this work, to prove hardness of refutations in the SoS proof system we show degree lower bounds of $d \in \Omega(n/\log n)$.
By the above relation, this implies exponential size lower bounds in $n$ as well. 
Similar size-degree trade-offs are also known for PC \cite{impagliazzo1999lower}, and thus it suffices for us to just show degree lower bounds to get exponential size lower bounds.
We denote the smallest maximum degree over all proofs that refute $\Axioms$ in PC and SoS, with $\Degree{Q \PC \bot}$ and $\Degree{Q \SOS \bot}$, respectively.
One motivation for proving lower bounds for algebraic proof systems, as opposed to propositional proof systems, is that often they imply lower bounds for a broad family of related algorithms for solving combinatorial optimisation problems.
Similarly, upper bounding the proof length has led to the fruitful discovery of many efficient algorithms.
The SoS proof system is of particular interest because of its close connection to the sum-of-squares hierarchy of semi-definite programming.
We refer the reader to the survey by Fleming, Kothari, and Pitassi \cite{fleming2019semialgebraic} for more details about the connections between the semi-algebraic proof systems and combinatorial optimisation.\par
In this work, we study the complexity of refuting \emph{perfect matchings} in PC and SoS.
Apart from being a natural problem in its own right, perfect matchings are also related to the pigeon hole principle \citep{maciel2000new, pitassi1993exponential, raz2004resolution, razbarov2002pgp, razborov2003resolution} and Tseitin formula \citep{filmus2013towards, galesi2019bounded, glinskih2017satisfiable, grigoriev2001linear}, two well studied formulae in proof complexity.
Assuming at most one pigeon fits in a single hole, the pigeon hole principle says $m$ pigeons cannot fit in $n < m$ holes.
If we construct the complete bipartite graph with the left vertices as $m$ pigeons and the right vertices as $n < m$ holes, proving the pigeon hole principle amounts to proving that such bipartite graph does not have a perfect matching.
There are other formulations of the pigeon hole principle (see the survey by Razborov \citep{razbarov2002pgp}), and almost all of them have short proofs in the sum of squares proof system.
In contrast, Tseitin formulae are known to require long proofs. The Tseitin formula over a graph claims that there is a spanning subgraph in which every vertex has odd degree.
If a graph has a perfect matching, then the subgraph described by the matching ensures that every vertex has odd degree.
However, formally refuting Tseitin formulae for expander graphs with an odd number of vertices, in the SoS proof system, requires degree linear in the number of vertices in the graph \cite{grigoriev2001linear}.
Given its close connections to the pigeon-hole and Tseitin, and the different behaviour of the two formulae, it is natural to determine the complexity of refuting perfect matchings for non-bipartite graphs.

To refute perfect matchings in an algebraic proof system, we first need to specify combinatorial constraints as algebraic equalities. Given an undirected graph $G=(V,E)$, $V = \{1, \ldots, n\}$, and a vector $\vec{b} = (b_1, \ldots, b_{n}) \in \Field^{n}$,
we define $\Card{G, \vec{b}}$ as the following set of polynomial constraints over variables $x_e$ for $e \in E$:

\[
        \Card{G, \vec{b}}=
        \Bigger[10]\{\begin{array}{@{}cl}
                x_{\highlight{e}}(1 - x_{\highlight{e}}) = 0 & \text{ for every $\highlight{e} \in E$}\\[3mm]
                \sum_{e \sim \highlight{v}} x_e = b_{\highlight{v}} & \text{ for every $\highlight{v} \in V$} \\[3mm]
        \end{array}
\]
For every $e \in E$, the equation $x_e(1 - x_e) = 0$ restricts the domain of the above variables to bits.
In plain words, $\Card{G, \vec{b}}$ denotes the claim that there exists a spanning subgraph $G' \subseteq G$ such that a vertex $v \in V(G)$ has $b_v$ edges incident to it in $G'$.
Note if there was an assignment of variables $(x_e)_{e \in E}$ that satisfies all equations in $\Card{G, \vec{1}}$, where $\vec{1} = (1, \ldots, 1) \in \Field^{n}$, it would imply that the graph $G$ has a perfect matching (given by the edges corresponding to variables with assignment 1).
Therefore, we define $\PM{G} := \Card{G, \vec{1}}$.
When $|V|$ is odd, $G$ trivially does not contain a perfect matching. How difficult is it to refute $\PM{G}$ in this case? In recent work, Austrin and Risse \cite{Austrin_2022} showed that refuting $\PM{G}$, in the Polynomial Calculus  and Sum-of-Squares system, in the case $G$ is a \emph{random $d$-regular graphs} with an odd number of vertices typically requires proofs with degree $\BigOmega{n/\log n}$.
They conjecture (see \citep[Section 6]{Austrin_2022}) that the hardness results should also apply to general expander graphs but leave showing so as an open problem.
In this work, we verify this by extending their result to all $d$-regular spectral expanders, that is, $d$-regular graphs with a mild condition on the spectral gap.
In fact, similar to Austrin and Risse, we reduce the hardness of refuting $\Card{G, \vec{t}}$, where $\vec{t} = (t, \ldots, t)$, for any odd value $t$, to the hardness of refuting $\Card{G, \vec{1}}$, where $\vec{1} = (1, \ldots, 1)$.
As another special case, this answers the \emph{even-colouring} case when $t = d/2$ is odd, a problem posed by Buss and Nordstr\"{o}m \citep[Open Problem 7.7]{buss2021proof}, which asks, \textit{``Are even colouring formulas over expander graphs hard for polynomial calculus over fields of characteristic distinct from 2 ?''}
Formally, we prove the following (for the definition of $(n, d, \lambda)$-graphs see Section \ref{sec:graph-theory-prelims}).

\begin{theorem}[Hardness Result For $\Card{G, \vec{t}}$]\label{thm:general-hardness-result}

  There exist universal constants $\epsilon, n_0, d_0 \in \Naturals$ such that for any odd $n \ge n_0$ and even $d \in [d_0, n]$, the following holds for \emph{any} $(n, d, \lambda)$-graph $G$ with $\lambda < \epsilon d$, and for any odd $1 \leq t \leq d$:  
  \[ \Degree{\Card{G, \vec{t}} \PC \bot} = \BigOmega{\nicefrac{n}{\log n}}\]
  \[ \Degree{\Card{G, \vec{t}} \SOS \bot} = \BigOmega{\nicefrac{n}{\log n}}\]  
\end{theorem}

\paragraph{Proof overview.}
We follow the overall approach of Austrin and Risse \cite{Austrin_2022}. Very briefly, the strategy is to obtain an affine restriction (see Definition \ref{def:affine-restriction}) $\Card{G, \vec{t}}|_\rho \equiv \PM{H}$ where $H$ is some graph for which refuting $\PM{H}$ requires large degree.
An example of such $H$ is given by Buss, Grigoriev, and Impagliazzo \cite{buss1999linear}. We now describe how to find such a restriction in more details: Using a result of Dragani\'c, Krivelevich, and Nenadov \cite{draganic22rolling}, we show that $H$ topologically embeds into a given expander graph $G$ with $\lambda \le \epsilon d$ for some universal small constant $\epsilon \in (0,1)$, such that all paths corresponding to the embedding have odd length.
The main technical ingredient of Austrin and Risse is also a similar embedding theorem, albeit a significantly more complicated one.
Moreover, we show that one can find such an embedding so that the subgraph of $G$ induced by vertices which are not part of the embedding has a perfect matching.
This allows us to use the restriction argument to transfer the hardness of $\PM{H}$ into the hardness of $\PM{G}$.
To extend this to hardness of $\Card{G, \vec{t}}$ for an odd $3 \le t \leq d$, it suffices to show that the graph $G'$ obtained from $G$ by removing all edges that participate in the embedding and the matching contains a $(t-1)$-regular spanning subgraph.
Austrin and Risse achieve this using the contiguity property of random regular graphs (and hence their hardness result for $\Card{G, \vec{t}}$ critically relies on randomness). 
Instead, we provide a significantly simpler and shorter argument based on Tutte's criterion.
 As a random $d$-regular graph is with high probability an $(n, d, \lambda)$-graph with $\lambda = \Theta(\sqrt{d})$ (see \cite[Theorem A]{tikhomirov2016spectralgapdenserandom}), our embedding theorem readily applies in the context of \citep{Austrin_2022}.\par
The rest of the document is structured as follows. 
In section \ref{sec:prelims}, we describe the requisite background from graph theory and proof complexity.
In Section \ref{sec:embed-machinery}, we describe the machinery for finding a desired topological embedding.
In Section \ref{sec:matching-machinery}, we prove conditions under which the residual graph has a perfect matching or, more generally, a $(t-1)$-regular spaning subgraph.
In Section \ref{sec:main-proof}, we use the tools from the previous sections to prove Theorem \ref{thm:general-hardness-result}. 
In Section \ref{sec:related-work} we briefly discuss a few other lower bounds using embeddings in proof complexity, and conclude with some future directions.

\section{Preliminaries}
\label{sec:prelims}

\subsection{Proof Complexity Preliminaries}
\label{sec:proof-system-prelims}

Let $\Axioms = \{ p_1 = 0, \dots, p_m = 0\}$ be a set of polynomial equations\footnote{The sum of squares proof system is a semi-algebraic proof system where $\Axioms$ may also contain inequalities of the form $p_i(x) \ge 0$. 
However, we only need equality constraints to express the existence of Perfect Matchings over graphs. 
Therefore to simplify our exposition, we write all our definitions using equality constraints only.}, which we refer to as axioms, over variables $\vec{X} = \{x_1, \dots, x_n, \bar{x}_1, \dots, \bar{x}_n\}$.

\begin{definition}[Sum Of Squares Refutation]\label{def:sum-of-squares} 
Given a set of $m$ polynomial equality constraints $\Axioms$ over the reals, a Sum of Squares (SoS) refutation is a sequence of polynomials $\Proof = (t_1, \dots, t_m; s_1, \dots, s_a)$ such that

\[ h \Def \sum_{i \in [m]} t_ip_i+ \sum_{i \in [a]} s_i^2 = -1\]

The degree of a proof $\Proof$ is

\[ \Degree{\Proof} \Def \max\left\{\max_{i \in [m]} \Degree{t_i} + \Degree{p_i}, \max_{i \in [a]} 2\Degree{s_i}\right\}\]

\end{definition}

Note that $s_i^2(x) \ge 0$ for any $x$ by definition. 
Therefore, if there were to exist some $x^{\star}$ such that $p_i(x^{\star}) = 0$ for all $p_i \in \Axioms$, then $\sum_{i \in [m]} t_i(x^\star)p_i(x^\star) = 0$. 
This would imply that $h \ge 0$, but our proof $\Proof$ shows that $h = -1$. 
Therefore, the existence of the set of polynomials $\Proof$ act as a formal proof of the claim that the set of polynomial equations in $\Axioms$ is unsatisfiable.
\begin{definition}[Complexity Of SoS Refutation]
If we let $\Pi$ denote the set of all valid SoS refutations for $\Axioms$, then the complexity of refuting $\Axioms$ in the SoS proof system is given by 
\[ \Degree{\Axioms \SOS \bot} \Def \min_{\Proof \in \Pi}\Degree{\Proof}\]

\end{definition}

Polynomial Calculus (PC) is a dynamic version of the static Nullstellensatz proof system \citep[Section 1.3]{fleming2019semialgebraic} operating over an arbitrary but fixed field, based on the following inference rules.
\begin{enumerate}
	\item From polynomial equations $f=0$ and $g=0$ where $f,g \in \Field[\vec{X}]$ we can derive $\alpha f + \beta g = 0$ for $\alpha, \beta \in \Field$.
	\item From polynomial $f=0$ where $f \in \Field[\vec{X}]$, we can derive $xf=0$ where $x \in \vec{X}$.
\end{enumerate}

\begin{definition}[Polynomial Calculus Refutations]\label{def:poly-calc-refutations}
A Polynomial Calculus (PC) refutation of $\Axioms$ over $\Field$ is a sequence of polynomials $\Proof = (t_1, \dots, t_l)$	such that $t_l = 1$, and for each $i \neq l$, either (1) $t_i \in \Axioms$, or (2) $t_i$ is derived from $\{t_j\}_{j < i}$ using the above rules.
The degree of the proof is given by $\Degree{\Proof} = \max_{i \in l}\Degree{t_i}$. If we let $\Pi$ denote the set of all PC refutations of $\Axioms$, then

\[ \Degree{\Axioms \PC \bot} \Def \min_{\Proof \in \Pi}\Degree{\Proof}\]
\end{definition}

To ensure Boolean variables, we assume the axioms $\Axioms$ always contain the equations $x_i^2 - x_i=0$ and $\bar{x}^2_i - \bar{x}_i = 0$ for all $i \in [n]$. 
Equivalently, we can also just work in the ring $\Field[x_1, \dots, \bar{x}_n]/(x_1^2 - x_1, \dots, \bar{x}_n^2 - \bar{x}_n)$ of multi-linear polynomials.
Multi-linearity implies that the degree of any proof can be at most $n$ i.e a proof of degree $\BigOmega{n}$ is the largest lower bound one can hope to achieve.
Additionally, we will also assume that $1 - x_i - \bar{x}_i=0$ is also included in $\Axioms$, for all $i \in [n]$, which ensures that the bar elements are bit complements of the non-bar elements.
The following lemma is by Buss, Grigoriev, and Impagliazzo \cite{buss1999linear} and gives an instance where perfect matching is hard to refute in the worst case. 

\begin{lemma}[Worst Case Hard Instance For PC]\label{lemma:worst-case-instance-PC}Given any odd $n \in \Naturals$, there exists a graph $H$ with $n$ vertices and maximum degree $\MaxDegree{H}= 5$ such that Polynomial Calculus over any field of characteristic different from 2 requires degree $\Theta(n)$ to refute $\Card{H, \vec{1}}$.
\end{lemma}
A description of the  worst case hard instance for SoS can be found in \citep[Theorem A.3]{Austrin_2022}.

\begin{lemma}[Worst Case Hard Instance For SOS]\label{lemma:worst-case-instance-sos}
Given any odd $n \in \Naturals$, there exists a graph $H$ with $n$ vertices and maximum degree $\MaxDegree{H}= 5$ such that SoS refutations requires degree $\Theta(n)$ to refute $\Card{H, \vec{1}}$.
\end{lemma}

\subsubsection{Affine restriction} 

An important lemma we will need is that given a set of axioms $\Axioms$ over the ring $\Field[x_1, \dots, x_n]$, a partial assignment of variables can only make refuting $\Axioms$ easier.
Given a set of $m$ polynomial equality constraints $\Axioms$ over boolean variables $\{x_1, \dots, x_n\}$, let the family of functions $\{f_i: \bit^n \rightarrow \{\True,\False\} \}_{i \in [m]}$, denote predicates for satisfiability for each constraint.
For example, given $\alpha \in \bit^n$, $f_i(\alpha) = \True$ if the $i$'th polynomial constraint $\axiom_i \in \Axioms$ is satisfied i.e $q_i(\alpha) = 0$.
We say $\Axioms$ is satisfied if $\exists \alpha \in \bit^n$ such that $f_i(\alpha) = \True \iff \axiom_i(\alpha)=0$ for all $i \in [m]$.
Given a map $\rho: \{x_1, \dots, x_n \} \rightarrow \{x_1, \dots, x_n, \bar{x}_1, \dots, \bar{x}_n, 1, 0 \}$, the restriction of a function $f: \bit^n \rightarrow \bit$, denoted by $f|_\rho$, is defined as $f|_\rho(x_1, \dots, x_n) = f(\rho(x_1), \dots, \rho(x_n))$.
Similarly, the restriction of formula $\Axioms$ is defined as $\Axioms|_\rho = \{f_{1}|_\rho, \dots, f_{m}|_\rho \}$.
Two formula $\Axioms$ and $\Axioms'$ are equivalent if they are element-wise equal, ignoring any functions that are constantly $\True$.
For example, $\Axioms = \{f_a, f_b, \True \}$ and $\Axioms' = \{f_a, f_b\}$ are equivalent, denoted as $\Axioms \equiv \Axioms'$.

\begin{definition}[Affine Restriction]\label{def:affine-restriction}
We say that an axiom $\Axioms'$ is an
affine restriction of $\Axioms$ if there is a map $\rho : \{x_1,\dots,x_n\} \rightarrow \{x_1, \dots, x_n, \bar{x}_1, \dots, \bar{x}_n, 1, 0 \}$ such that $\Axioms \equiv \Axioms|_{\rho}$.
\end{definition}

\begin{lemma}\label{lemma:affine_restriction}
Let $\Axioms, \Axioms'$ be axioms such that $\Axioms'$ is an affine restriction of $\Axioms$, and each axiom
of $\Axioms$ depends on a constant number of variables, then
\begin{enumerate}
	\item For any field $\Field$ it holds that $\Degree{\Axioms \PC \bot} \in \BigOmega{\Degree{\Axioms' \PC \bot}}$
	\item $\Degree{\Axioms \SOS \bot} \in \BigOmega{\Degree{\Axioms' \SOS \bot}}$
\end{enumerate}
\end{lemma}

The proof for the above lemma can be found in \citep[Lemma 2.2]{Austrin_2022}.
What the above lemma says is that if we have a graph $G$ with odd vertices with constant degree, that has a perfect matching on a subset of even vertices on the graph, then the size of the proof to refute $\PM{G}$ is at least as large as refuting a perfect matching in $G$ with the even vertices removed.

\subsection{Graph Theory Preliminaries}
\label{sec:graph-theory-prelims}

We use standard graph theoretic notation. For a graph $\Graph$, we use $\Vertices{\Graph}$ and $\Edges{\Graph}$ to denote the vertices and edges of $\Graph$. For a vertex $v \in \Vertices{\Graph}$, we use $\Neighbourhood{G}{v} = \{ u \in \Vertices{\Graph} : (u,v) \in E' \}$ to denote the neighbourhood of $v$ in $G$, and $\degree{G}{v} := |\Neighbourhood{G}{v}|$. Given two sets $S, T \subseteq \Vertices{G}$, we  use $\CutEdges{S}{T}{G}$ to denote the number of edges in $G$ with one endpoint in $S$ and one endpoint in $T$.  Note that we do not require $S$ and $T$ to be disjoint; in case they are not disjoint, every edge with both endpoints in $S \cap T$ is counted twice in $\CutEdges{S}{T}{G}$. If the graph $G$ is clear from the context, we omit the subscript. Given $W \subseteq V(G)$, we denote with $G[W]$ the subgraph of $G$ induced by $W$. We say that a subgraph $G' \subseteq G$ is \emph{spanning} if $V(G') = V(G)$.

Next, we give a definition of pseudorandom graphs.

\begin{definition}[$\EnDeeLambda$-graphs]\label{def:expander-graphs}
Let $G$ be a $d$-regular graph on $n$ vertices, and, let $\lambda_1 \geq \lambda_2 \ge \dots \geq \lambda_n$ denote eigenvalues of the adjacency matrix of $G$.
We say $G$ is an $\EnDeeLambda$-graph if $\ExpansionFactor{\Graph} := \underset{{\{2, \dots, n\}}}{\max}|\lambda_i| \leq \lambda$.
\end{definition}

The following is a well known result of Alon and Chung \cite{alon88mixing}.

\begin{lemma}[Expander Mixing Lemma]\label{lemma:expanders-mixing-lemma}
  Given an $\EnDeeLambda$-graph $G$, for any $S, T \subseteq \Vertices{G}$ we have
$$
  \left| \CutEdges{S}{T}{G} - \frac{d}{n}|S||T| \right|
  \le \lambda\sqrt{\Size{S}\Size{T} }
$$
\end{lemma}

We make use of the following two well known criteria of Tutte \cite{tutte1952factors, tutte1947factorization}. 
Note that both of the lemmata ask for properties which are stronger than what Tutte criteria requires\footnote{More specifically, in the Tutte criterion $q(G)$ denotes the number of \emph{odd} sized connected components.}; however, they are easier to state and verify in our application. We denote with $q(G)$ the number of connected components in a graph $G$.

\begin{lemma}[Tutte's Criterion]\label{lemma:tutte-criterion}
If a graph $G$ has \emph{even} number of vertices and for every subset $S \subseteq \Vertices{G}$ we have $\OddComponents{G \setminus S} \leq |S|$, then $G$ contains a perfect matching.
\end{lemma}

\begin{lemma}[Tutte's Generalised Criterion]\label{lemma:tutte-criterion-factor}
  Let $f \in \Naturals$ be even. Suppose $G$ is a graph such that for every pair of disjoint sets $S, T \subseteq V(G)$ the following holds:

  \[ \OddComponents{G \setminus (S \cup T)} \le |S|f - \sum_{w \in T} (f - \Size{\Neighbourhood{G}{w} \setminus S}) \]
  Then $G$ contains a spanning subgraph $G' \subseteq G$ which if $f$-regular.
\end{lemma}


\subsection{Probabilistic Tools}

\begin{lemma}[Multiplicative Chernoff bound]\label{lemma:mult-chernoff}
Suppose $X_1, ..., X_n$ are identical independent random variables taking values in $\{0, 1\}$. Let $X$ denote their sum and let $\mu = n\Mean{}{X_1}$ denote the sum's expected value. Then, for any $0 < \delta < 1$, we have

\[ \Prob{|X - \mu| \geq \delta \mu} \leq 2\exp(-\delta^2\mu/3)\]

\end{lemma}

A dependency graph for a set of events $E_1, . . . , E_n$ is a graph $G=(V, E)$ such that $V = \{1,.. . , n\}$ and,  for $i= 1,\dots, n$, event $E_i$ is mutually independent
of the events $\{E_j | (i, j) \notin E\}$. The degree of the dependency graph is the maximum degree of any vertex in the graph.

\begin{lemma}[Lov\`asz Local Lemma]\label{lemma:lll}
Let $E_1,...,E_n$ be a set of events, and assume that the following hold:
\begin{enumerate}
\item The degree of the dependency graph given by $(E_1, \dots, E_n)$ is bounded by $d$.

\item For all $i \in [n]$, $\PProb{E_i}{\Dist} \leq \beta.$

\item $\beta \leq \frac{1}{4d}$.
\end{enumerate}
Then
\[ \Pr{\overset{n}{ \underset{i=1}{\bigcap}} \hspace{0.1cm}  \overline{E_i}}{\Dist} > 0\]

\end{lemma}




The proof of \nameref{lemma:mult-chernoff} and \nameref{lemma:lll}
can be found in any textbook on randomised algorithms (for example, see \citep[Chapter 1, Chapter 7]{mitzenmacher2017probability}). 
The following lemma is by \citep[Lemma 4.3]{Austrin_2022}, re-derived here for completeness.

\begin{lemma} \label{thm:partition}
For every $0 < c < 1$ and $\gamma> 0$, there exists $d_0$ such that the following holds. If $G$ is a $d$-regular graph, for some $d \ge d_0$, then there exists a subset $A \subseteq V(G)$ such that
$$
    cd - \gamma d \leq \left| \Neighbourhood{G}{v} \cap A \right| \leq cd + \gamma d 
$$
for every $v \in V(G)$.
\end{lemma}
\begin{proof}
We prove the existence of such a partition $A \subseteq V(G)$ using the probabilistic method.
For each $v \in \Vertices{\Graph}$, we toss an independent coin $X_i$ with bias $c$.
We include $v$ in $A$ if and only if $X_i = 1$. Thus, $\vec{X} \Def (X_1, \dots, X_n) \in \bit^n$ is a random variable that describes how we choose $A$. For any $v \in V$, let $Y_v \Def |\Neighbourhood{G}{v} \cap \RedA|$ denote the random variable that counts the number of neighbours of $v$ in $A$.
Define $\delta \Def \gamma / c$, and for every $v \in V$ let $E_v = \Indicator{|Y_v - d c | \geq \delta c d}$ denote the bad event that $v$ has too many or too few neighbours in $A$.
Observe that the dependency graph of events $\{ E_v \}_{v \in V}$ has maximum degree at most $d^2$ (only vertices at most two hops away from $v$ affect how many of $v$'s neighbours are in $A$; there at most $d^2$ such vertices).
As $\Graph$ is $d$-regular, $\Mean{\vec{X}}{Y_v} = c d$. By the \nameref{lemma:mult-chernoff}, for any $v \in V$ we have $\PProb{E_v}{\vec{X}} \leq 2\exp(-\delta^2c d/3) =: \beta$.
For $d$ sufficiently large we have $\beta \leq 1/(4d^2)$, and so $\beta d^2 \leq 1/4$. All the conditions of \nameref{lemma:lll} are satisfied, from which we conclude that, with positive probability, none of the bad events happen. This implies the desired $A \subseteq V(G)$ exists.
\end{proof}




\section{Topological embedding}
\label{sec:embed-machinery}

In this section we describe the topological embedding result of Dragani\'c, Krivelevich, and Nenadov \cite{draganic22rolling}. We start with a necessary definition.

\begin{definition}[Sub-divisions]\label{def:subdivisions}
Given a graph $H$ and a function $\sigma: \Edges{H} \rightarrow \Naturals$, the $\sigma$-subdivision of $H$, denoted by $\Subdivision{H}{\sigma}$, is the graph obtained by replacing each edge in $\Edges{H}$ with a path of length $\sigma(e)$ joining the end points of $e$ such that all these paths are mutually vertex disjoint, except at the end points.
\end{definition}

If a graph $G$ contains $H^{\sigma}$ for some $\sigma: \Edges{H} \to \Naturals$, then we say $G$ contains $H$ as a \emph{topological minor}. In our application, it will be important that we can control the parity of $\sigma(e)$. The following result follows directly from \cite[Theorem 1]{draganic22rolling}.

\begin{theorem}[Embedding Theorem] \label{thm:embedding_top}
For every $D \in \mathbb{N}$ there exist $\alpha, \xi, C > 0$, such that the following holds. Suppose $G$ is a graph with $n$ vertices and $m \ge Cn$ edges such that for every pair of disjoint subsets $S, T \subseteq V(G)$ of size $|S|, |T| \ge \xi n$, we have
\begin{align*}
  \bigl| e_{G}(S, T) - |S||T|p \bigr| \le \xi |S||T| p \label{eq:props}
\end{align*}
where $p = m / \binom{n}{2}$. Then $G$ contains $H^{\sigma}$, where $H$ is any graph with maximum degree at most $D$, $H^\sigma$ has at most $\alpha n$ vertices, and $\sigma(e) \ge \log n$ for every $e \in E(H)$.
\end{theorem}

When $G$ is an $\EnDeeLambda$ graph, we will we make use of Theorem \ref{thm:embedding_top} to show that $G$ satisfies the required properties, thereby contains $H$ as a topological minor.
This gives us the following corollary.
\begin{corollary} \label{cor:embedding_top}
  For every $D \in \mathbb{N}$ there exist $d_0, n_0 \in \Naturals, \epsilon, \alpha \in (0,1)$, such that the following holds. 
  Suppose $G$ is an $(n,d,\lambda)$-graph where $d \ge d_0$, and $\lambda < \epsilon d$, and $n \ge n_0$. 
  Let $B \subseteq V(G)$ be a subset of size $|B| \ge n/20$, and $H$ is any graph with maximum degree at most $D$ and at most $\alpha \frac{n}{\log n}$ vertices.
Then the induced sub-graph $G[B]$ contains $H^{\sigma}$ such that $\sigma(e)$ is odd for every $e \in E(H)$.
\end{corollary}

\begin{proof}
  Let $m$ denote the number of edges in the induced subgraph $G[B]$, which gives us $2m = \CutEdges{B}{B}{G}$.
Denote $b \Def \Size{B}$ and Define $p = m / \binom{b}{2}$. By the \nameref{lemma:expanders-mixing-lemma}, we have

\[ \left| 2m - \frac{d}{n}b^2 \right| \le \lambda b \]
Dividing both sides with $b(b-1)$, and observing that $\frac{db^2}{n(b-1)b} = \frac{d}{n}\frac{b}{b-1} = \frac{d}{n}(1 + \frac{1}{b-1})$, we further get

\begin{align}
\left| \frac{2m}{b(b-1)} - \frac{d}{n} - \frac{d}{n(b-1)} \right| &\le \frac{\lambda}{b-1}  \label{eq:eml-post-process}
\end{align}

From this we have 
\begin{align}
  \left| p - \frac{d}{n}\right| &= \left| p - \frac{d}{n} - \frac{d}{n(b-1)} + \frac{d}{n(b-1)}\right|  \\
                                &\le \left| p - \frac{d}{n} - \frac{d}{n(b-1)} \right| + \frac{1}{(b-1)} \label{eq:triangle}\\
                                &\le \frac{\lambda}{(b-1)} + \frac{1}{(b-1)}  \label{eq:final} \\
  & \le \frac{2\lambda}{b}
\end{align}

Let us briefly justify each step: \eqref{eq:triangle} comes from the triangle inquality and  $d/n \le 1$;
Equation \eqref{eq:final} comes from Equation \eqref{eq:eml-post-process};
the last inequality comes from the assumption that $b \in \BigOmega{n}$, so the inequality holds for $n$ large enough.\\

Let $\xi$ be as given by the \nameref{thm:embedding_top}. Using the bound on the difference between $p$ and $d/n$, for every disjoint subsets $S, T \subseteq B$ of size $|S|,|T| \ge \xi n$, for $\lambda < \epsilon d$ where $\epsilon$ is sufficiently small, we have

\begin{align*}
  \Bigl| \CutEdges{S}{T}{G}  - p|S||T| \Bigr| &\le \left| \CutEdges{S}{T}{G} - \frac{d}{n}|S||T| \right| + \left| \frac{d}{n}|S||T| - p|S||T| \right| \\
                                    &\le \lambda \sqrt{|S||T|} + \frac{2 \lambda}{b} |S||T| \le \xi |S||T|p 
    \end{align*}

With the lower bounds on $S,T$ and $B$, we can make $\epsilon$ sufficiently small with respect to $\xi$ to get the upper bound in the last step.
    Let $\sigma: E(H) \to \Naturals$ be the constant function where $\sigma(e)$ is the smallest odd integer larger than $\log n$. 
As $G[B]$ has at least $Cn$ edges (by the \nameref{lemma:expanders-mixing-lemma}), $\sigma(e) \le 2 + \log n$, and $H$ has at most $\alpha \frac{n}{\log n}$ vertices, we can invoke the \nameref{thm:embedding_top} to conclude that $G[B]$ contains $H^\sigma$.
\end{proof}

\section{Perfect matching and regular subgraphs}
\label{sec:matching-machinery}
As described earlier, the second ingredient in our hardness proof is showing that a certain residual graph contains a perfect matching or a spanning $(t-1)$-regular subgraph. In this section we state and prove these ingredients.

\begin{lemma}[Perfect Matching Lemma]\label{thm:perfect-matching}
Let $G$ be an  $(n, d, \lambda)$-graph with $\lambda < d/50$, and suppose $G' \subseteq G$ satisfies $\delta(G') \ge 0.9d$. Then for all $S \subseteq V(G')$, $G' \setminus S$ has at most $|S|$ connected components, that is, $q(G' \setminus S) \le |S|$. Therefore, if $G'$ has an even number of vertices then it contains a perfect matching.
\end{lemma}

\begin{proof}
  Let $U = V(G')$.  We aim to show that the graph $G' \setminus S$ has at most $|S|$ connected components. If $|S| \ge |U|/2$ then $G' \setminus S$ has at most $|S|$ vertices, so the the upper bound on connected components trivially holds. For the remainder of the proof we can assume $|S| < |U|/2$.

  \noindent
  We claim the following:
\begin{equation}
  \text{For every partition } X \cup Y = U \setminus S,\text{ with } |X|,|Y| \ge |S|/3,
  \quad
  \CutEdges{X}{Y}{G'} \ge 1 \implies  q(G' \setminus S) \le |S|.
\label{eq:mini-claim}
\end{equation}
\flushleft
To see why, assume towards a contradiction that there exists an edge in $G'$ between \emph{every} partition $X \cup Y = U \setminus S$, where $|X|, |Y| \ge \Size{S}/3$, \emph{and} $G' \setminus S$ has more than $|S|$ connected components. Denote the  vertex sets of these components by $C_1, \ldots, C_k$, for some $k > |S|$.
Let $X^\star := C_1 \cup \ldots \cup C_{s}$ and $Y^\star := C_{s+1} \cup \ldots \cup C_{k}$, where $s = \lfloor |S|/2 \rfloor \ge |S|/3$.
By construction, even if each component $C_i$ is a singleton set, we get that  $\Size{X^\star}, \Size{Y^\star} \ge \Size{S}/3$.
Now as all $C_i$'s are disjoint connected components, there can be no edge between $X^\star$ and $Y^\star$. Therefore, we have found a partition $X^\star \cup Y^\star = U \setminus S$ with $|X^\star|, |Y^\star| \ge |S|/3$ without an edge between them, which contradicts our assumption that all appropriately sized partitions have at least one edge between them.

To complete our proof, it suffices to show $e_{G'}(X, Y) \geq 1$ for every partition $X \cup Y = U \setminus S$ with $|X|, |Y| \ge |S|/3$.\\
\flushleft
Consider some arbitrary partition $X \cup Y$ of $U \setminus S$, with $|X|, |Y| \ge |S|/3$, and without loss of generality assume $|X| \le |Y|$. 
Then by a simple counting argument we get
\begin{equation}
|X| \le (|U| - |S|)/2 \le (n - |S|)/2  \label{eq:upp-X}
\end{equation}

We have: 
\begin{align}
\CutEdges{X}{X}{G'} + \CutEdges{X}{S}{G'} &\le \CutEdges{X}{X}{G} + \CutEdges{X}{S}{G} \nonumber \\
   &\le \frac{d}{n}|X|^2 + \lambda |X| + \frac{d}{n}|X||S| + \lambda\sqrt{|X||S|}        \\
&\le \frac{d}{n} |X| \frac{(n-|S|)}{2} + \lambda |X| + \frac{d}{n}|X||S| + \lambda\sqrt{3} |X| \label{eq:upper-bound-X}  \\
&< \frac{d|X|}{2} + \frac{d|X||S|}{2n} + 3 \lambda |X|\\
&< \frac{d}{2}|X| + \frac{d}{4}|X| + 3 \lambda |X| \label{eq:smaller-than-one} \\
&< \cTilde|X|  \label{eq:upper-bound-S}
\end{align}
These steps are justified as follows: The first equation follows from the \nameref{lemma:expanders-mixing-lemma}; Equation \eqref{eq:upper-bound-X} comes from Equation \eqref{eq:upp-X} and $|X| \ge |S|/3$;
Equation \eqref{eq:smaller-than-one} comes from $\Size{S} < \frac{n}{2}$; and  Equation \eqref{eq:upper-bound-S} comes the assumption $\lambda < d/50$.
By the assumption $\delta(G') > 0.9d$ we conclude that there is an edge in $G'$ with one vertex in $X$ and the other in $V(G') \setminus (X \cup S) = Y$.
\end{proof}



The next lemma shows that subgraphs of $(n, d, \lambda)$-graphs with large minimum degree contain regular spanning subgraphs.

\begin{lemma}[Regular Subgraph Lemma]\label{lemma:f-factor}
  For every $C > 1$ there exists $d_0 = d_0(C)$ such that the following holds.
Suppose $G$ is an $\EnDeeLambda$ graph with $\lambda < \epsilon d$ and $d \ge d_0$, where $\epsilon < 1/100C^{3/2}$.
If $G' \subseteq G$ has minimum degree $\MinimalDegree{G'} \ge d - C$, then $G'$ contains a spanning $f$-regular subgraph for any even $2 \le f \le d/2$. 
\end{lemma}

\begin{proof}
  We prove this lemma using the \nameref{lemma:tutte-criterion-factor}.
  We need to show that for any pair of disjoint sets $S, T \subseteq \Vertices{G'}$, we have

  \begin{align}
    \OddComponents{G' \setminus (S \cup T)} &\le \Size{S}f - \sum_{w \in T} (f - \Size{\Neighbourhood{G'}{w} \setminus S}) \label{eq:main-goal}
  \end{align}

  As $\epsilon < \frac{1}{100C^{3/2}}$ and $C > 1$, we have that $\epsilon < 1/100$. 
This implies that $G$ is an $\EnDeeLambda$ graph with $\lambda < d/100$.
We set $d_0 := d_0(C)$ large enough such that for all $d \ge d_0$, even after deleting at most $C$ edges incident on each vertex of the $d$-regular graph $G$ to get $G'$, we have the minimum degree of $G'$ to be $\MinimalDegree{G'} \ge  d - C > 9d/10$.
Therefore the conditions of Lemma \ref{thm:perfect-matching} are satisfied,
thus
$$\OddComponents{G' \setminus (S \cup T)} \le \Size{S \cup T} = |S| + |T|$$
To prove \eqref{eq:main-goal}, it suffices to show 
  \begin{align}
    |S| + |T|  &\le \Size{S}f - \Size{T}f + \Size{T}(d-C) - \CutEdges{S}{T}{G'}  \label{eq:to-show}  \\ 
    &\le |S|f - \sum_{w \in T} (f - \Size{\Neighbourhood{G'}{w} \setminus S}) \nonumber
  \end{align}  
We distinguish a few cases.

  \paragraph{Case 1:} Suppose $\Size{S} \leq \Size{T}$. As $f \leq d/2$, we have 
  \begin{align}
    \Size{S}f + \Size{T}(d-C-f)  &\ge (\Size{S} + \Size{T})(d/2 -C) \label{eq:to-show-1}
  \end{align}

  The condition described by the inequality \eqref{eq:to-show} is satisfied via the following analysis.
\begin{align}
  |S| + |T| + \CutEdges{S}{T}{G'} &\le  |S| + |T| + |S||T|\frac{d}{n} + \epsilon d \sqrt{\Size{S}\Size{T}} \label{eq:eml}\\
                                 &\le |S| + |T| + \frac{d}{4}(\Size{S} + \Size{T}) + \epsilon d \frac{1}{2}(\Size{S} + \Size{T}) \label{eq:algebra}\\
                                 &\le (\Size{S} + \Size{T})\left(\frac{d}{4} + 1 + \frac{
                                   \epsilon d}{2}\right) \\
                                 &<  (\Size{S} + \Size{T})\left(d/2-C\right) \label{eq:eps-assumption}\\
                                 &\le \Size{S}f - \Size{T}(d-C-f) \label{eq:plug}
 \end{align}  
  
Equation \eqref{eq:eml} comes from the \nameref{lemma:expanders-mixing-lemma} and $\lambda < \epsilon d$, together with an obvious upper bound $e_{G'}(S, T) \le e_G(S, T)$.
Equation \eqref{eq:algebra} comes from the fact that $\Size{S}\Size{T} \leq \left(\frac{\Size{S}+\Size{T}}{2}\right)^2 \le n\frac{(\Size{S}+\Size{T})}{4}$.
Equation \eqref{eq:eps-assumption} comes from $\epsilon < \nicefrac{1}{100C^{3/2}}$, $C >1$ and $d_0$ being sufficiently large. Equation \eqref{eq:plug} follows from Equation \eqref{eq:to-show-1} which gives us what we want.

\paragraph{Case 2:} Suppose $\Size{S} > \Size{T}$. As $f \ge 2$, we have
  \begin{align}
    \Size{S}f +  \Size{T}(d-C-f)  &\ge 2\Size{S} +\Size{T}(d-C-2) \label{eq:to-show-2}
  \end{align}
  
  To show Equation \eqref{eq:to-show}, it suffices to show that
\begin{align}
  \CutEdges{S}{T}{G'} \leq \Size{S} + \Size{T}(d-C-3) \label{eq:to-show-3}
\end{align}

Now we distinguish between two subcases.

\begin{enumerate}

\item{
    If $\Size{T} \le \frac{\Size{S}}{C+3}$, then Equation \eqref{eq:to-show-3} follows from a trivial bound $\CutEdges{S}{T}{G'} \leq \Size{T}d$.
  }

\item{
$ \frac{\Size{S}}{C+3} < \Size{T} < \Size{S}$. As $\Size{S} + \Size{T} < n$ we have $\Size{S} < n - \frac{\Size{S}}{C+3}$, thus $|S| < n\left(\frac{C+3}{C+4}\right)$.
Using the \nameref{lemma:expanders-mixing-lemma}, we have
\begin{align}
\CutEdges{S}{T}{G} \le \CutEdges{S}{T}{G'} &\le \frac{d}{n}\Size{S}\Size{T} + \epsilon d \sqrt{\Size{S}\Size{T}} \nonumber	\\[8pt]
&< d\frac{(C+3)}{(C+4)}\Size{T} + \epsilon d \Size{T}\sqrt{(C+3)} \nonumber  \\
&= |T|d \left( \frac{C+3}{C+4} + \varepsilon \sqrt{C+3} \right) \nonumber \\
&< |T|d\left( \frac{C+3}{C+4} + \frac{1}{2(C+4)} \right) = |T|d\left( 1 - \frac{1}{2(C+4)}\right) \nonumber
\end{align}
where the penultimate inequality follows from the upper bound on $\varepsilon$. For $d$ sufficiently large in terms of $C$ we obviously have
$$
    |T|d \left(1 - \frac{1}{2(C+4)}\right) < |T|d - |T|(C + 3) < |S| + |T|(d - C - 3)
$$
hence \label{eq:to-show-3} is satisfied.


  }
  
\end{enumerate}
\end{proof}

\section{Proof of Theorem \ref{thm:general-hardness-result}}
\label{sec:main-proof}

In this section we prove Theorem \ref{thm:general-hardness-result}. As $\Card{G, \vec{t}} \equiv \Card{G, \vec{d-t}}$, without loss of generality we only prove the theorem for $t \le d/2$. 

\begin{proof}[Proof of Theorem \ref{thm:general-hardness-result}]
Let $G=(V,E)$ be an $\EnDeeLambda$-graph on an odd number of vertices with $\lambda < \epsilon d$, where $\epsilon < \nicefrac{1}{100C^{3/2}}$ and $C=6$. 
For sufficiently small constant $\alpha \in (0,1)$, let $\HardInstance$ denote the graph on $h=\alpha\frac{n}{\log n}$ vertices as given by Lemma \ref{lemma:worst-case-instance-sos} (to show lower bounds for PC, we use $H$ from Lemma \ref{lemma:worst-case-instance-PC}). Recall that any SoS proof which refutes $\PM{H}$ has degree $\BigOmega{h}$. We now make use of $H$ to show the hardness of refuting $\Card{G, \vec{t}}$. The idea is to find a restriction $\rho$ such that $\Card{G, \vec{t}}|_\rho \equiv \PM{H}$. We achieve this through the following steps.

\begin{enumerate}
\item{Invoke Lemma  \ref{thm:partition} with parameters $c=0.925$ and $\gamma = 0.025$ to get subsets $A \subseteq V(G)$ and $B = \Vertices{G} \setminus A$, such that for every $u \in V(G)$ we have
    
\begin{align} 
0.9d  &\leq   \Size{\Neighbourhood{G}{u} \cap \RedA} \leq 0.95d \label{eq:size_nbr_A}\\[5pt]
0.05d  &\leq   \Size{\Neighbourhood{G}{u} \cap \GreenB} \leq 0.1d 
\label{eq:size_nbr_B}
\end{align}
}

\item{From equations \eqref{eq:size_nbr_A} and \eqref{eq:size_nbr_B}, $\Size{B} > n/20$, with room to spare and $\Size{E(G[B])} \ge nd/20$.
    By Corollary \ref{cor:embedding_top}, $G[B]$ contains $H^{\sigma}$ such that each $\sigma(e)$ is odd. Let us denote a subgraph of $G[B]$ corresponding to $H^{\sigma}$ by $G_{\EmbeddingFunc}$.
    We can describe $G_\EmbeddingFunc$ as a function $\EmbeddingFunc: V(H) \to B$ together with a collection of pairwise vertex-disjoint (other than at the endpoints) paths  $\Path{\Embedding{u}}{\Embedding{v}}$ in $G[B]$, for $(u,v) \in E(H)$.
    Observe that it is at least as hard to refute\footnote{Note that this by itself does not guarantee it is hard to refute $\PM{G}$. We need item (3) and this to show hardness of refuting $\PM{G}$.} $\PM{G_{\EmbeddingFunc}}$ as it is to refute $\PM{H}$.
To see why, let $y_1, \dots, y_{E(H)}$ denote the variables for the $\PM{H}$ formulae for each edge of $H$.
We use as shorthand $\mathcal{Y} = (y_e)_{e \in E(H)}$ and $\Complement{\mathcal{Y}} = (\Complement{y}_e)_{e \in E(H)}$.
Define a mapping $\rho': \Edges{G_\EmbeddingFunc} \to \{0,1, \mathcal{Y}, \Complement{\mathcal{Y}} \}$ as follows. For each $(u,v) \in E(H)$, let $\rho'(x_e) = y_{u,v}$ where $e$ is the first edge on the path $\Path{\Embedding{u}}{\Embedding{v}}$.
Map each variable $x_e$ for $e \in \Path{\Embedding{u}}{\Embedding{v}}$ alternately to $y_{u,v}$ or $\bar{y}_{u,v}$, such that the edges of the path adjacent to $\Embedding{u}$ and $\Embedding{v}$ are set to $y_{u,v}$. This is always possible as $\Path{u}{v}$ has odd length. Observe that $\PM{G_\EmbeddingFunc}|_{\rho'} \equiv \PM{H}$. 
  }

\item{ As $n$ is odd and $\Size{\Vertices{G_{\EmbeddingFunc}}}$ is odd, we have that $U = \Vertices{G} \setminus \Vertices{G_{\EmbeddingFunc}}$ has even size. From \eqref{eq:size_nbr_A} we have that $G[U]$ has minimum degree at least $\cTilde$. As $\lambda < d/50$ (with room to spare), we can invoke Lemma \ref{thm:perfect-matching} to conclude $G[U]$ has a perfect matching $M$. 
}
  
\item{
  Consider the subgraph $G' \subseteq G$ obtained by deleting all edges $e \in \Edges{G_{\EmbeddingFunc}} \cup M$, where $M$ is the perfect matching from the step above. As $\MaxDegree{H} \leq 5$, every vertex $u \in G$ loses at most $5$ edges in this process. Thus, we have $\MinimalDegree{G'} \ge d - 5$. As $\lambda < \epsilon d$, by the \nameref{lemma:f-factor} we have that $G'$ contains a $(t-1)$-regular spanning subgraph $G''$. }

\end{enumerate}

We finally define $\rho$ as follows:
\[
\rho(x_e) =
\left\{
\begin{array}{@{}cl}
\rho'(e) & \text{if } e \in \Edges{G_{\EmbeddingFunc}} \\
1 & \text{if } e \in M \cup E(G'') \\
0 & \text{otherwise.}
\end{array}
\right.
\]
Then $\Card{G, \vec{t}}|_\rho \equiv \PM{H}$, thus our theorem follows from Lemma \ref{lemma:affine_restriction}.
\end{proof}

\section{Related Work And Concluding Remarks}
\label{sec:related-work}

In propositional\footnote{As opposed to algebraic proof complexity} proof complexity, there a few prior examples of the strategy of embedding a worst case instance into a host graph to show lower bounds for a larger class of objects \cite{itsykson2021Near, pitassi2016PolyLogFrege}.
Pitassi, Rossman, Servedio and Tan show Tseitin lower bounds for Frege proof systems \citep{pitassi2016PolyLogFrege} by relying on the embedding result by Rubinfeld and Kleinberg \citep{kleinberg1996short}, which allows one to embed any bounded degree graph $H$ of size $\BigO{n/\poly(\log n)}$ into an expander graph on $n$ vertices as a minor (not necessarily a topological one). 
Krivelevich and Nenadov simplify and improve the above embedding theorem to allow for embedding any graph $H$ with size $\BigO{n/\log n}$ as an ordinary minor \citep{krivelevich2021completeMinors}.
However, embedding a hard instance $H$ into $G$ as an ordinary minor does not guarantee that the hardness of $H$ is preserved in the setting considered in this paper. In particular, it is entirely possible that one of the edge contractions to obtain the minor results in $H$ now being easy to refute.
Thus, these embedding theorems cannot be directly applied.
Instead, as described in section \ref{sec:main-proof}, one way to preserve hardness is to use  embedding theorems that allow for topological embeddings that allow for edge sub-divisions of odd size \citep{draganic22rolling, nenadov2023routing}.
In order to get a topological embedding, Austrin and Risse modify the ordinary embedding theorem in \citep{krivelevich2021completeMinors} but critically rely on the host graph being random.
In this work, we use the embedding theorem by Dragani\'c, Krivelevich, and Nenadov \cite{draganic22rolling}, which greatly simplifies the argument. Moreover, we avoid the use of the contiguity argument present in \cite{Austrin_2022} by directly utilising Tutte's criterion and the Expander Mixing Lemma.

\par
In summary, we show degree lower bounds for refuting $\Card{G, \vec{t}}$ for odd $t$ in $\EnDeeLambda$ graphs in the SoS and PC proof systems.
There is still a $\log n$ gap between the largest possible proof in such systems, and our lower bounds (similar to Austrin and Risse).
It is not inherently clear that such a gap should exist.
The gap is an artefact of $d$ being constant, which makes the graphs sparse i.e we need $\Theta(\log n)$ edges to form a path between any two nodes.
This implies, that $\BigOmega{n/\log n}$ is the largest hard instance we can topologically embed in any graph.
Thus, if the worst case lower for refuting perfect matchings was indeed $\Omega(n)$, we would need a more direct proof of the statement without using a smaller hard instance.
We leave the issue of resolving the tightness of our lower bound as an open problem for future work.

\bibliographystyle{abbrv}
\bibliography{perfect-matching}
\clearpage
\appendix

\end{document}